\documentclass[a4paper,11 pt,reqno]{amsart}
\usepackage{a4,amsbsy,amscd,amssymb,amstext,amsmath,latexsym,oldgerm,enumerate,verbatim}

\makeatletter                                       
\def\l@subsection{\@tocline{2}{0pt}{2.5pc}{2.5pc}{}}
\makeatother                                        

\makeatletter                                                  
\def\chapter{\clearpage\thispagestyle{plain}\global\@topnum\z@ 
\@afterindenttrue \secdef\@chapter\@schapter}
\makeatother

\newtheorem{thmgl} {Theorem}    
\newtheorem{propgl}{Proposition}
\newtheorem{lemgl} {Lemma}

\newtheorem{lemnn}{Lemma}

\newtheorem{cornn}{Corollary}

\theoremstyle{definition}

\newtheorem{remgl} {Remark}
\newtheorem{remsgl} [remgl]{Remarks}

\newcommand{\mf}{\mathfrak}
\newcommand{\mc}{\mathcal}
\newcommand{\mb}{\mathbb}

\newcommand{\ov}{\overline}

\newcommand{\sm}{\setminus}         

\newcommand{\ot}{\otimes}           
\newcommand{\la}{\langle}
\newcommand{\ra}{\rangle}

\newcommand{\Mat}{{\rm Mat}}
\newcommand{\Der}{{\rm Der}}

\newcommand{\Ker}{{\rm Ker}}

\newcommand{\rank}{{\rm rank}}

\newcommand{\tr}{{\rm tr}}
\newcommand{\id}{{\rm id}}
\newcommand{\gr}{{\rm gr}}   

\newcommand{\g}{\mf{g}}
\newcommand{\h}{\mf{h}}
\let\ttie\t
\newcommand{\tie}[1]{{\let\t\ttie \ttie#1}}
\renewcommand{\t}{\mf{t}}  

\renewcommand{\b}{\mf{b}}
\renewcommand{\c}{\mf{c}}

\newcommand{\gl}{\mf{gl}}
\newcommand{\spl}{\mf{sl}}

\newcommand{\Ad}{{\rm Ad}}              
\newcommand{\Lie}{{\rm Lie}}
\newcommand{\GL}{{\rm GL}}
\newcommand{\PGL}{{\rm PGL}}
\newcommand{\SL}{{\rm SL}}

\newcommand{\SU}{{\rm SU}}

\newcommand{\SO}{{\rm SO}}
\newcommand{\PSO}{{\rm PSO}}
\newcommand{\Spin}{{\rm Spin}}
\newcommand{\Sp}{{\rm Sp}} 
\newcommand{\PSp}{{\rm PSp}}

\newcommand{\Pic}{{\rm Pic}}

\begin{document}

\title{Infinitesimal invariants in a function algebra}
\author{Rudolf Tange}

\begin{abstract}
Let $G$ be a reductive connected linear algebraic group over an algebraically closed field of positive characteristic and let $\g$ be its Lie algebra. First we extend a well-known result about the Picard group of a semisimple group to reductive groups. Then we prove that, if the derived group is simply connected and $\g$ satisfies a mild condition, the algebra $K[G]^\g$ of regular functions on $G$ that are invariant under the action of $\g$ derived from the conjugation action, is a unique factorisation domain.
\end{abstract}

\address{School of Mathematics, University of Southampton, Highfield, SO17 1BJ, UK}
\email{rtange@maths.soton.ac.uk}
\maketitle

\section*{Introduction}
The study of noncommutative algebras that are finite over their centres has raised the interest in those centres and commutative algebras that look like those centres. Important examples are the universal enveloping algebra of a reductive Lie algebra $\g$ in positive characteristic and quantised enveloping algebras at a root of unity. Their centres are not smooth, but have several nice properties. For example, they are normal domains and admit a "Friedlander-Parshall factorisation" into two well understood algebras (see e.g. \cite{BrGo} Thm~3.5 (5) and Thm.~4.1 (3)). Analogues in a commutative setting for these centres are the invariant algebras $K[\g]^\g$ and $K[G]^\g$ (see below for precise definitions). In these cases the Friedlander-Parshall factorisation was obtained by Friedlander and Parshall (\cite{FrPa}) and Donkin (\cite{D}, one can replace $G_1$ by $\g$).

In this note we will consider the question whether such a commutative algebra is a unique factorisation domain (UFD). For the algebra $K[G]^\g$ this question was posed by S.\ Donkin. In the noncommutative setting of the universal enveloping algebra this question was answered in \cite{PrT} for type $A_n$. This question turned out to be rather trivial for the algebra $K[\g]^\g$ (see \cite{PrT} 3.3, 3.4). For completeness we state some general results in the final section. To answer this question for $K[G]^\g$ we need to know when $K[G]$ is a UFD. This is well-known if $G$ is semi-simple simply connected. In Section~\ref{s.picard} we extend this result to reductive groups. Our main result is that $K[G]^\g$ is a UFD if the derived group is simply connected and $\g$ satisfies a mild condition. This is proved in Section~\ref{s.liealginvs}. Another question one could ask is whether such an algebra has a rational (i.e. purely transcendental over the ground field) field of fractions. Even in the commutative setting this question is unsolved outside type $A_n$. In the final section we state an affirmative answer to this question for the algebra $K[\SL_n]^{\spl_n}$.

Throughout $K$ is an algebraically closed field of characteristic $p>0$, $G$ is a connected reductive algebraic group over $K$ and $\g=\Lie(G)$ is its Lie algebra. The conjugation action of $G$ on itself induces an action of $G$ on $K[G]$, the algebra of regular functions on $G$.
We will refer to this action and its derived $\g$-action as {\it conjugation} actions. We also use this terminology for the induced actions of subgroups of $G$ resp. subalgebras of $\g$ on localisations of $K[G]$. The conjugation action of $G$ on $K[G]$ is by algebra automorphisms, so the conjugation action of $\g$ on $K[G]$ is by derivations. In the notation of \cite{Bo} I.3.7 it is given by $x\mapsto*x-x*:\g\to\Der_K(K[G])$: the left invariant vector field determined by $x$ minus the right invariant vector field determined by $x$.

If $L$ is a Lie algebra over $K$ and $V$ is an $L$-module, then $v\in V$ is called an $L$-{\it invariant} if $x\cdot v=0$ for all $x\in L$ and an $L$-{\it semi-invariant} if the subspace $Kv$ is $L$-stable. If $v\in V$ is a nonzero $L$-semi-invariant, then there is a unique linear function $\lambda:L\to K$ such that $x\cdot v=\lambda(x)v$ for all $x\in L$ and for this $\lambda$ we have $\lambda([L,L])=0$.

\section{The Picard group of $G$}\label{s.picard}

Fix a Borel subgroup $B$ of $G$ and a maximal torus $T$ of $G$ that is contained in $B$. We denote the character group of $T$ by $X(T)$. The usual pairing of characters and cocharacters of $T$ is denoted by  $\la\,,\,\ra$. Let $B^-$ be the opposite Borel subgroup relative to $T$ and let $U$ and $U^-$ be the unipotent radicals of $B$ and $B^-$. By a "root" we will mean a root (i.e. nonzero weight) of $T$ in $\g$. The root subgroup associated to a root $\alpha$ is denoted by $U_\alpha$. We denote $K$ as an additive group by ${\mb G}_a$. If $\theta_\alpha:{\mb G}_a\stackrel{\sim}{\to}U_\alpha$ is an isomorphism of algebraic groups, then
\begin{equation}
t \theta_\alpha(a)t^{-1}=\theta_\alpha(\alpha(t)a)\text{\qquad for all }t\in T\text{ and }a\in K.
\label{eq.T-on-rootsubgroup}
\end{equation}
A root $\alpha$ is called positive ($\alpha>0$) if it is a root of $T$ in $\Lie(B)$. Let $\alpha_1,\ldots,\alpha_s$ be the simple roots and let $W=N_G(T)/T$ be the Weyl group of $G$ relative to $T$. To every root $\alpha$ there is associated an element $s_\alpha\in W$ that is a reflection on $\mb R\ot_{\mb Z}X(T)$. For $w\in W$ define $U^-_w:=wU^-w^{-1}\cap U^-$. Clearly $wU^-_{w^{-1}}w^{-1}=U^-_w$. In the terminology of \cite{Bo}~14.3 (see also 14.12) we have that the unipotent group $U^-_w$ is directly spanned by the $U_{-\alpha}$ with $\alpha>0$  and $w^{-1}(\alpha)>0$ in any order. By a version of the Bruhat decomposition (\cite{Bo} 14.12) we have that $G=\bigcup_{w\in W}B^-wB$ disjoint and $B^-wB=U^-_wwB$ for every $w\in W$. If $\dot w$ is a representative in $N_G(T)$ of $w$, then the map
\begin{equation}
(u^-,t,u)\mapsto u^-\dot w tu:U^-_w\times T\times U\to B^-wB
\label{eq.celliso}
\end{equation}
is an isomorphism of varieties. This implies that the codimension of $BwB$ in $G$ is $l(w)$, the length of $w$ with respect to our choice of simple roots. Furthermore we have that
\begin{equation}
\dot wU_\alpha\dot w^{-1}=U_{w(\alpha)}.
\label{eq.W-on-rootsubgroup}
\end{equation}
Denote the "big cell" $B^-B$ by $\Omega$ and denote the closure of $B^-s_{\alpha_i}B$ by $\Gamma_i$. As is well known, $\Omega$ is open and the $\Gamma_i$ are the irreducible components of $G\sm\Omega$.

As $G$ is a smooth variety, we can and will identify $\Pic(G)$ with the divisor class group ${\rm Cl}(G)$ of $G$. For a nonzero rational function $f\in K(G)$ we denote, as usual, the divisor (of zeros and poles) of $f$ on $G$ by $(f)$. For $\chi\in X(T)$ we define the function $\ov\chi\in K[\Omega]\subseteq K(G)$ by
$$\ov\chi(u^-tu)=\chi(t)$$
for $u^-\in U^-$, $t\in T$ and $u\in U$.

The next theorem corrects Theorems~2 and 6 in \cite{Po} and extends \cite{FoIv} Cor.~4.5 (and \cite{FoIv} Cor~4.6 and \cite{Iv} Prop.~3.5) to reductive groups. See also \cite{Mar} Cor~3 to Prop.~4 (there $G$ is assumed to be semi-simple simply connected).

It is not clear how the method of the proof of Theorem~2 in \cite{Po} can be extended to reductive groups. The point is that to associate to certain Weil divisors on $G$ a projective representation of $G$, one needs to know that the invertible regular functions on $G$ are the nonzero constant functions. For reductive groups that are not semi-simple this is no longer true. We will use a method which is more elementary and much shorter. Besides the proof of (iii), of which the main assertion (iv) is an easy consequence, all arguments are standard (see e.g. \cite{Che} and \cite{Po}) and are only given for completeness.

Recall that $T\cap DG$ is a maximal torus of the derived group $DG$ and that its character group is naturally identified with a sublattice of the weight lattice of the root system of $DG$ relative to $T\cap DG$ (which is the same as that of $G$ relative to $T$).

\begin{thmgl}[see also \cite{FoIv}, \cite{Iv}, \cite{Po}]\label{thm.picard}
Let $L$ be the group of divisors that has the $\Gamma_i$ as a basis, let $P$ be the weight lattice of the root system of $G$ relative to $T$ and let $\varpi_i$ be the $i^{\rm th}$ fundamental weight. Then the following holds:
\begin{enumerate}[{\rm(i)}]
\item $L$ consists of the divisors that are fixed by $B$ under right multiplication and by $B^-$ under left multiplication.
\item For $f\in K(G)$ nonzero we have $(f)\in L$ if and only if $f$ is regular and nowhere zero on $\Omega$ if and only if $f$ is a scalar multiple of $\ov\chi$ for some $\chi\in X(T)$.
\item For $\chi\in X(T)$ we have $(\ov\chi)=\sum_{i=1}^s\la\chi,\alpha_i^\vee\ra\Gamma_i$.
\item The homomorphism $P\to \Pic(G)$ that maps $\varpi_i$ to the divisor class of $\Gamma_i$ is surjective and induces an isomorphism $P/X(T\cap DG)\stackrel{\sim}{\to}\Pic(G)$.
\end{enumerate}
\end{thmgl}

\begin{proof}
(i).\ Let $D$ be a divisor that is fixed by $B$ under right multiplication and by $B^-$ under left multiplication. Then the same holds for the support $Supp(D)$ of $D$. If $Supp(D)$ would intersect $\Omega$, then it would contain all of $\Omega$ which is impossible. Therefore $Supp(D)\subseteq G\sm\Omega$ which means that $D\in L$.\\
(ii).\ We have $(f)\in L$ if and only if $Supp((f))$ does not intersect $\Omega$ if and only if $f$ is regular and nowhere zero on $\Omega$ if and only if $f$ is an invertible element of $K[\Omega]$. Now $K[\Omega]\cong K[T]\ot_KK[U^-\times U]$ which is a polynomial algebra over $K[T]$. The invertible elements of this algebra are the invertible elements of $K[T]$ which are, as is well-known and easy to see, the scalar multiples of the characters of $T$. The assertion now follows from the fact that the embedding $K[T]\hookrightarrow K[\Omega]$ corresponding to the embedding $K[T]\hookrightarrow K[T]\ot_KK[U^-\times U]$ is given on $X(T)$ by $\chi\mapsto\ov\chi$.\\
(iii).\ Since $\ov\chi$ is regular and nowhere zero on $\Omega$ we have that $(\ov\chi)=\sum_{i=1}^sn_i\Gamma_i$ for certain integers $n_i$. To determine the $n_i$ we follow the proof of Prop.~II.2.6 in \cite{Jan}. Choose isomorphisms $\theta_\alpha:{\mb G}_a\stackrel{\sim}{\to}U_\alpha$ as in \cite{Jan} II.1.1-1.3 (there denoted by $x_\alpha$). See also \cite{St2} \S 3, Lemma~19 and Cor.~6 to Thm.~$4'$. For $a\in K^\times$ define, as in \cite{Jan}  and \cite{St2}, $n_\alpha(a):=\theta_\alpha(a)\theta_{-\alpha}(-a^{-1})\theta_\alpha(a)$. Then $n_\alpha(a)=\alpha^\vee(a)n_\alpha(1)$ is a representative of $s_\alpha$ in $N_G(T)$ and
\begin{equation}
n_\alpha(1)=\alpha^\vee(a^{-1})\theta_\alpha(a)\theta_{-\alpha}(-a^{-1})\theta_\alpha(a)
\text{\qquad for }a\in K^\times\text{ and }\alpha\text{ a root.}
\label{eq.nalpha}
\end{equation}
Let $i\in\{1,\ldots,s\}$. Since $\Gamma_i$ intersects the open set $\Omega_i:=s_{\alpha_i}\Omega=U^-_{s_{\alpha_i}}s_{\alpha_i}U_{-\alpha_i}B$, it suffices to show that the restriction $(\ov\chi)\vert_{\Omega_i}$ of $(\ov\chi)$ to $\Omega_i$ equals $\la\chi,\alpha_i^\vee\ra\Gamma_i\cap\Omega_i$.
Take $n_{-\alpha_i}(1)$ as a representative in $N_G(T)$ of $s_{\alpha_i}$. From \eqref{eq.W-on-rootsubgroup} it follows that conjugation with $n_{-\alpha_i}(1)$ induces an automorphism of $U^-_{s_{\alpha_i}}$. From this and \eqref{eq.celliso} with $\dot w=1$ we deduce that
\begin{equation}
(u^-,a,t,u)\mapsto u^-n_{-\alpha_i}(1)\theta_{-\alpha_i}(a)tu:U^-_{s_{\alpha_i}}\times\mb{G}_a\times T\times U\to\Omega_i
\label{eq.specialcelliso}
\end{equation}
is an isomorphism of varieties. Let $u^-\in U^-_{s_{\alpha_i}}$, $t\in T$, $u\in U$ and $a\in K^\times$. Replacing $\alpha$ and $a$ in \eqref{eq.nalpha} by $-\alpha_i$ and $-a$ we obtain
$$
n_{-\alpha_i}(1)\theta_{-\alpha_i}(a)=
\alpha_i^\vee(-a)\theta_{-\alpha_i}(-a)\theta_{\alpha_i}(a^{-1})=
\theta_{-\alpha_i}(b)\alpha_i^\vee(-a)\theta_{\alpha_i}(a^{-1})
$$
for some $b\in K$. Here we used that $(-\alpha_i)^\vee(c)=\alpha_i^\vee(c^{-1})$ for $c\in K^\times$.
So
\begin{align*}
&\ov\chi(u^-n_{-\alpha_i}(1)\theta_{-\alpha_i}(a)tu)=
\ov\chi(u^-\theta_{-\alpha_i}(b)\alpha_i^\vee(-a)t\theta_{\alpha_i}(\alpha_i(t)^{-1}a^{-1})u)=\\
&\chi(\alpha_i^\vee(-a)t)=
(-a)^{\la\chi,\alpha_i^\vee\ra}\chi(t)
\end{align*}
Thus \eqref{eq.specialcelliso} gives us an isomorphism of $\Omega_i$ with an affine space with some coordinate hyperplanes removed and $\ov\chi$, considered as a rational function on $\Omega_i$, is equal to an invertible regular function on $\Omega_i$ times the $\la\chi,\alpha_i^\vee\ra^{\rm th}$ power of minus the coordinate functional that defines $\Gamma_i\cap\Omega_i$. So $(\ov\chi)\vert_{\Omega_i}=\la\chi,\alpha_i^\vee\ra\Gamma_i\cap\Omega_i$.\\
(iv).\ Denote the homomorphism by $F$. Since $\Omega$ is affine and $K[\Omega]$ is isomorphic to a localised polynomial algebra, we have $\Pic(\Omega)={\rm Cl}(\Omega)=0$. So for every divisor $D$ on $G$ there exists a nonzero $f\in K(G)$ such that $(f)\vert_\Omega=D\vert_\Omega$, that is, such that $D-(f)\in L$. So every divisor $D$ on $G$ is linearly equivalent to a divisor in $L$. This means that $F$ is surjective. By (ii) and (iii) the kernel of $F$ is the image of $\chi\mapsto\sum_{i=1}^s\la\chi,\alpha_i^\vee\ra\varpi_i:X(T)\to P$, which is $X(T\cap DG)$ in its usual identification with a sublattice of $P$.
\end{proof}

\begin{cornn}
$K[G]$ is a unique factorisation domain if and only if $\Pic(G)=0$ if and only if the derived group of $G$ is simply connected.
\end{cornn}

\begin{remsgl}\label{rem.picard}
1.\ The fact that $n_\alpha(a)$ is a representative of $s_\alpha$ in $N_G(T)$ implies that $s_{\alpha}T(U_\alpha\sm\{e\})\subseteq TU_\alpha U_{-\alpha}$. From this it follows that $s_{\alpha_i}\Gamma_j=\Gamma_j$ for $i\ne j$ and therefore that $(s_{\alpha_i}\Omega)\cap\Gamma_j=\emptyset$ if $i\ne j$, as suggested by the above proof.\\
2.\  As in \cite{Che} Prop.~4, it follows now that for any divisor $D$ and any $g\in G$ the divisors $D$, $gD$ and $Dg$ are linearly equivalent: It is clear that for any divisor $D$ and any $g\in B$ we have $Dg\sim D$, since any divisor is linearly equivalent to one in $L$ and the divisors in $L$ are fixed by $B$ under right multiplication. Now let $g\in G$. Then there exists an $h\in G$ such that $hgh^{-1}\in B$. So we have that $Dh^{-1}\sim Dh^{-1}(hgh^{-1})=Dgh^{-1}$ and therefore that $D\sim Dg$. The proof that $D$ and $gD$ are linearly equivalent is completely analogous.\\
3.\ The mistake in Theorem~2 in \cite{Po} is caused by the fact that the varieties $\Delta_i$ in \cite{Po} are not the same as those in \cite{Che}. In \cite{Po} $\Delta_i$ is the closure of $B^-s_{\alpha_i}B/B$. In \cite{Che} $\Delta_i$ is the closure of $Bs_{\alpha_i}w_0B/B=Bw_0s_{\alpha_i^*}B/B=w_0B^-s_{\alpha_i^*}B/B$. Making the corresponding replacements one obtains a correct proof of \cite[Theorem~2]{Po}.

That Theorem~6 in \cite{Po} is incorrect can easily be seen by applying it to $\GL_n$. This gives $\Pic(\GL_n)\cong\mb{Z}_n$ if $p\nmid n$, which is clearly not true. The isogeny $\pi:~\tilde\GL_n~\to~\GL_n$ is inseparable if $p\mid n$ and then it does not induce an isomorphism $\tilde\GL_n/\Ker(\pi)\cong\GL_n$. Furthermore we have that $\Ker(\pi)\nsubseteq D\tilde\GL_n$ if $n$ is not a power of $p$.\\
4.\ It might be possible to prove Theorem~1(iv) by examining the exact sequence of \cite{FoIv} Prop.~4.2 (or \cite{Iv} Prop.~2.6) and taking $G'$ the cover of $G$ with $\Pic(G')=0$ that is constructed in these papers.
\end{remsgl}

\section{Infinitesimal invariants in $K[G]$}\label{s.liealginvs}

To proof the main result we have to impose a mild condition on $\g$. We will first investigate this condition. Recall that the rank of $G$ is the dimension of a maximal torus of $G$. For any $x\in\g$ we have $\dim(C_G(x))\geq\rank(G)$. An element $x$ of $\g$ is called {\it regular} if its $G$-orbit under the adjoint action has dimension $\dim(G)-\rank(G)$ or, equivalently, if its centraliser $C_G(x)$ has dimension $\rank(G)$. It is well known that the set of regular elements is open and nonempty and that the set of regular semi-simple elements is open. Since the set of semi-simple elements is constructible, it follows that $\g$ has regular semi-simple elements if and only if the semi-simple elements are dense in $\g$.
Let $T$ be a maximal torus of $G$. It is easy to see that $x\in\Lie(T)$ is regular if and only if the differentials of the roots of $G$ relative to $T$ are nonzero at $x$. From this and the fact that every semi-simple element of $\g$ is tangent to a (maximal) torus (\cite{Bo}~11.8), it follows that $\g$ has no regular semi-simple elements if and only if for some (and therefore every) maximal torus the differential of some root relative to $T$ is zero.

The next lemma shows that the existence of regular semi-simple elements in $\g$ is a mild assumption. Note that $\Sp_2\cong\SL_2$.
\begin{lemgl}\label{lem.regsselt}
If $\g$ has no regular semi-simple elements, then $p=2$ and the derived group $DG$ has a quasi-simple factor that is isomorphic to $\Sp_{2n}$ for some $n\ge1$. Conversely, if $p=2$ and $G$ contains some $\Sp_{2n}$, $n\ge1$, as a direct algebraic group factor, then $\g$ has no regular semi-simple elements.
\end{lemgl}

\begin{proof}
Let $T$ be a maximal torus of $G$ and let $\Phi$ be the root system of $G$ relative to $T$. By assumption some $\alpha\in\Phi$ has differential zero. Let $H$ be the quasi-simple factor of $DG$ whose root system $\Psi$ is the irreducible component of $\Phi$ containing $\alpha$. Then $\alpha=p\chi$ for some $\chi\in X(T\cap H)$, so $2=\la p\chi,\alpha^\vee\ra=p\la\chi,\alpha^\vee\ra$ and therefore $p=2$. Now choose a basis $(\alpha_1,\ldots,\alpha_r)$ of $\Psi$ such that $\alpha_1=\alpha$. Then we have $\la\alpha_1,\alpha_j^\vee\ra=2\la\chi,\alpha_j^\vee\ra$ for all $j\in\{1,\ldots,r\}$. So some row of the Cartan matrix (with our basis the first) has even entries. Inspecting the tables in \cite{Bou} we see that $\Psi$ has to be of type $A_1$ or of type $C_n$,  $n\ge2$. If $H$ would be of adjoint type then the above basis of $\Psi$ would also be a basis of the character group of $T\cap H$ and the differential of $\alpha$ would be nonzero. Since the connection index is $2$ this means that $H$ is simply connected and therefore isomorphic to $\Sp_{2n}$ for some $n\ge1$.

If $H\cong\Sp_{2n}$ for some $n\ge1$, then twice the first fundamental weight is the highest root. So the second statement follows from the fact that an equality $\alpha=2\chi$ on $T\cap H$, for some direct algebraic group factor $H$ of $G$ and some root $\alpha$ of $H$, is also valid on $T$ if we extend $\chi$ to $T$ in the obvious way.
\end{proof}

The condition in the first statement of the preceding lemma is not sufficient for the non-existence of semi-simple elements in $\g$, since $\gl_2$ (in fact any $\gl_n$) always has regular semi-simple elements. If $G=\SL_2$ and $p=2$, then the semi-simple elements of $\g$ are the multiples of the identity. If $G=\Sp_{2n}$ ($n\ge2$) and $p=2$, then the derived algebra $[\g,\g]$ is equal to the sum of the Lie algebra of a maximal torus and the root spaces corresponding to the short roots (see \cite{St1} sect.~2 at the end). Since $[\g,\g]$ is stable under the adjoint action we then have that the set of semi-simple elements of $\g$ is contained in the proper subspace $[\g,\g]$ of $\g$.

Note that the proposition below implies that every $\g$-semi-invariant of $\g$ in $K[G]$ is a $\g$-invariant.

\begin{propgl}\label{prop.stableideal}
Assume that $\g$ has regular semi-simple elements. Let $f\in K[G]$ be a  regular function. If the ideal $K[G]f$ is stable under the conjugation action of $\g$ on $K[G]$, then $f$ is $\g$-invariant.
\end{propgl}

\begin{proof}
Note that the assumption that the ideal $K[G]f$ is $\g$-stable is equivalent to the assumption that $f$ divides $x\cdot f$ for all $x\in\g$. So there exists a linear function $F:\g\to K[G]$ such that
\begin{equation}
x\cdot f=F(x)f\text{\quad for all }x\in\g
\label{eq.F}
\end{equation}
Clearly we may assume that $f\ne0$. We will show that $F=0$.

Let $T$ be a maximal torus of $G$ on which $f$ does not vanish and let  $\Phi$ be the root system of $G$ relative to $T$.
Let $U,U^-,U_\alpha,\Omega$ and the choice of the positive roots be as in Section~\ref{s.picard}. Note that $\Omega$ is stable under the conjugation action of $T$.  For $\alpha\in\Phi$ let $\theta_\alpha:{\mb G}_a\stackrel{\sim}{\to}U_\alpha$ be an isomorphism of algebraic groups. Choose an ordering $\beta_1,\ldots,\beta_N$ of the positive roots. Then the multiplication defines  isomorphisms of varieties $\prod_{i=1}^NU_{\beta_i}\stackrel{\sim}{\to}U$ and $\prod_{i=1}^NU_{-\beta_i}\stackrel{\sim}{\to}U^-$. Using the above isomorphisms and the isomorphism \eqref{eq.celliso} for $\dot w=1$, we obtain for each positive (resp. negative) root $\alpha$ a coordinate function $\xi_\alpha\in K[\Omega]$: first project to $U$ (resp. $U^-$), then project to $U_\alpha$ and then apply $\theta_\alpha^{-1}$. Composing the elements of a basis of the character group of $T$ with the projection $\Omega\to T$ we obtain $r=\rank(G)$ further coordinate functions $\zeta_1,\ldots,\zeta_r\in K[\Omega]$. The functions $\xi_\alpha$, $\alpha\in\Phi$, and $\zeta_i$, $i\in\{1,\ldots,r\}$, are together algebraically independent and $K[\Omega]=A[(\xi_\alpha)_{\alpha\in\Phi}]$, where $A=K[(\zeta_i)_i][(\zeta_i^{-1})_i]$, a localised polynomial algebra.

As a polynomial algebra over $A$, $K[\Omega]$ has a grading. Using \eqref{eq.T-on-rootsubgroup} it is easy to see that the grading of $K[\Omega]$ is stable under the conjugation action of $T$ and therefore also under the conjugation action of $\Lie(T)$. Now let $h\in\Lie(T)$. Taking degrees on both sides of \eqref{eq.F} we obtain that either $F(h)=0$ or $\deg(F(h))=0$. So $F(h)\in A$, which means that $F(h)(u^-tu)=F(h)(t)$ for all $u^-\in U^-,t\in T, u\in U$. Now $(h\cdot f)\vert_T=h\cdot (f\vert_T)=0$, since the conjugation actions of $T$ and $\Lie(T)$ on $K[T]$ are trivial. So \eqref{eq.F} gives $F(h)\vert_T=0$, since $f\vert_T\neq0$. Therefore $F(h)=0$.

Put $\mf S=\bigcup\{\Lie(T)\,\vert\,f\vert_T\neq0\}$, the union of the Lie algebras of the maximal tori of $G$ on which $f$ does not vanish. We have shown that $F$ is zero on $\mf S$. So it suffices to show that $\mf S$ is dense in $\g$ in the Zariski topology. Denote the set of zeros of $f$ in $G$ by $Z(f)$. Fix a maximal torus $T$ of $G$. Let $\mu:G\times T\to G$ be the morphism that maps $(g,t)$ to $gtg^{-1}$. Put $\mc O=\{g\in G\vert gTg^{-1}\nsubseteq Z(f)\}$. This set is nonempty, since $Z(f)$ is a proper closed subset of $G$ and the union of the conjugates of $T$ is the set of semi-simple elements of $G$ which contains an open dense subset of $G$ (see \cite{Bo} 11.10). Furthermore $\mc O= pr_G(\mu^{-1}(G\sm Z(f)))$ is an open set, because the projection $pr_G:G\times T\to G$ is open. Let $\nu:G\times\Lie(T)\to\g$ be the morphism that maps $(g,x)\mapsto\Ad(g)(x)$. Its image is the set of semi-simple elements of $\g$ which contains an open dense subset of $\g$. But then $\mf S=\nu(\mc O\times\Lie(T))$ also contains an open dense subset of $\g$, since it is a constructible set and $\nu$ is continuous.

\end{proof}

\begin{thmgl}\label{thm.ufd}
Assume that the derived group of $G$ is simply connected and that $\g$ has regular semi-simple elements. Then the invariant algebra $K[G]^\g$ is a unique factorisation domain. Its irreducible elements are the irreducible elements of $K[G]$ that are invariant under $\g$ and the $p$-th powers of the irreducible elements of $K[G]$ that are not invariant under $\g$.
\end{thmgl}

\begin{proof}
By Theorem~\ref{thm.picard} $K[G]$ is a UFD. Let $f$ be a nonzero element in $K[G]^\g$ and suppose $f=f_1f_2$ where $f_1,f_2\in K[G]$ are coprime and not a unit. Then we have for all $x\in\g$ that $(x\cdot f_1)f_2=-f_1(x\cdot f_2)$. Since $K[G]$ is a UFD and $f_1$ and $f_2$ are coprime this means that $f_i$ divides $x\cdot f_i$ for $i=1,2$ and all $x\in\g$. By Proposition~\ref{prop.stableideal} it follows that $f_1,f_2\in K[G]^\g$.

Now suppose $f=g^n$ for some $n\in \mb N$. Write $n=sp+r$ with $s,r\in {\mb Z}_+$ and $0\le r<p$. Then $0=x\cdot f=ng^{n-1}(x\cdot g)$. For $r\ne 0$ this yields $g\in K[G]^\g$, while for $r=0$ we have  $f=(g^p)^s$ and obviously $g^p\in K[G]^\g$.

This shows that any irreducible element in $K[G]^\g$ is either an irreducible element of $K[G]$ invariant under $\g$ or a $p$-th power of an irreducible element in $K[G]\setminus K[G]^\g$. Now the unique factorisation property of $K[G]^\g$ follows from that of $K[G]$.
\end{proof}

\begin{remgl}
Note that by Remark~\ref{rem.picard}.3 the units of $K[G]^\g$ are the same as those of $K[G]$, since characters of $K[G]$ are obviously invariant under the conjugation representation of $G$ and $\g$.
\end{remgl}

\section{Complements}

The following proposition is an obvious generalisation of  Proposition~1 in \cite{PrT}. Note that if $A_0=K$, the second assumption says that every $L$-semi-invariant is an $L$-invariant.

\begin{propgl}\label{prop.filtration}
Let $L$ be a Lie algebra over $K$ and let $A$ be a commutative $K$-algebra on which $L$ acts by derivations. Assume that $A$ is a unique factorisation domain and that it has an $L$-stable $K$-algebra filtration $A_0\subseteq A_1\subseteq A_2\cdots$ whose associated graded is an integral domain. Assume furthermore that every $a\in A$ such that $A_0a$ is $L$-stable is an $L$-invariant. Then the following holds:
\begin{enumerate}[{\rm (i)}]
\item If, for $a\in A$, the ideal $Aa$ is $L$-stable, then $a$ is an $L$-invariant.
\item The invariant algebra $A^L$ is a unique factorisation domain and the irreducible elements of $A^L$ are  the irreducible elements of $A$ that are invariant under $L$ and the $p$-th powers of the irreducible elements of $A$ that are not invariant under $L$.
\end{enumerate}
\end{propgl}

\begin{proof}
(i).\ We may assume that $a\ne 0$. We have $a$ divides $x\cdot a$ for all $x\in L$~(*). Since the associated graded is an integral domain, the (filtration) degree of a product is the sum of the degrees of the factors. Taking degrees in (*) we obtain that $L\cdot a\subseteq A_0a$ which means, by assumption, that $a$ is an $L$-invariant.\\
(ii).\ This is precisely as in the proof of Theorem~\ref{thm.ufd}.
\end{proof}

\begin{lemgl}\label{lem.semiinv}
Every $\g$-semi-invariant in $K[G]$ is a $\g$-invariant.
\end{lemgl}

\begin{proof}
We first make the following observation: Let $V$ be a restricted $\g$-module and let $v\in V$ be a $\g$-semi-invariant. If there exists a maximal torus $T$ of $G$ such that $\Lie(T)\cdot v=0$, then $v$ is a $\g$-invariant. This follows from the fact that an element of a root space of $\g$ relative to $T$ is nilpotent and therefore has to act nilpotently on any finite dimensional sub $\g$-module of $V$ (a nilpotent endomorphism of $Kv$ is zero).

Now let $f\in K[G]$ be a nonzero semi-invariant for $\g$ and let $\lambda$ be the linear functional on $\g$ such that $x\cdot f=\lambda(x)f$ for all $x\in\g$. As we saw at the end of the proof of Proposition~\ref{prop.stableideal}, there exists a maximal torus $T$ of $G$ such  that $f\vert_T\ne0$. Again as in the proof of Proposition~\ref{prop.stableideal}, we use the fact that the conjugation action of $\Lie(T)$ on $K[T]$ is trivial and obtain that $\lambda(x)=0$ for all $x\in\Lie(T)$.
\end{proof}

The first statement of the next proposition is a generalisation of \cite{PrT} Lemma~2. The condition of the second statement only fails if $p=2$ and the the derived group has a quasi-simple factor which is adjoint and has root system of type $A_1$ or $B_n$, $n\ge2$ (see the proof of Lemma~\ref{lem.regsselt}). The third statement includes one case not covered by Theorem~\ref{thm.ufd}: $\SL_2$ in characteristic $2$.

\begin{propgl}\label{prop.application}
Proposition~\ref{prop.filtration} with $A_0=K$ and $L=\g$ applies in the following three cases:
\begin{enumerate}[{\rm(1)}]
\item If $\g$ has regular semi-simple elements. For $A=K[\g]$ under the adjoint action.
\item If the differential of a coroot is never zero. For $A=S(\g)$ under the adjoint action.
\item If $G=\SL_n$ and $A=K[G]$.
\end{enumerate}
\end{propgl}

\begin{proof}
(1).\  As filtration we choose the filtration that comes from the natural grading of $K[\g]=S(\g^*)$. The associated graded is obviously isomorphic to $K[\g]$. It remains to prove that every $\g$-semi-invariant in $K[\g]$ is a $\g$-invariant. We use the observation in the proof of Lemma~\ref{lem.semiinv}. Let $f\in K[G]$ be a nonzero semi-invariant for $\g$. Let $T$ be a maximal torus of $G$ and let $\t$ be its Lie algebra. As we have seen at the beginning of Section~\ref{s.liealginvs}, the assumption on $\g$ is equivalent to the assumption that the set of semi-simple elements is dense in $\g$. But this set is the union of the conjugates of $\t$ (see \cite{Bo}~11.8). So, after replacing $T$ by a conjugate, we may assume that $f\vert_\t\ne0$. Using the fact that the adjoint action of $\t$ on $K[\t]$ is trivial we deduce, as in the proof of Lemma~\ref{lem.semiinv}, that $\t\cdot f=0$.\\
(2).\ Here we also use the filtration that comes from the natural grading. It remains to prove that every $\g$-semi-invariant in $S(\g)=K[\g^*]$ is a $\g$-invariant. Again we use the observation in the proof of Lemma~\ref{lem.semiinv}. Let $f\in K[G]$ be a nonzero semi-invariant for $\g$. Let $T$ be a maximal torus of $G$ and let $\t$ be its Lie algebra. Put $h_\alpha=d\alpha^\vee(1)\in\t$. By assumption the $h_\alpha$ are nonzero, so there exists a $\lambda\in\h^*$ that is nonzero at each $h_\alpha$. We extend $\lambda$ to a linear functional on $\g$ by requiring it to be zero on the root spaces. Denote the subspace of all linear functionals on $\g$ with this property by $\t'$. Note that $T$ fixes $\t'$ pointwise. Since $[\g_\alpha,\g_{-\alpha}]$ is spanned by $h_\alpha$ and $\c_\g(\lambda)$ is $T$-stable, it follows that $\c_\g(\lambda)=\t$ and that $C_G(\lambda)^0=T$.
By the same arguments as \cite{Bo}~11.10 we have that the set ${\rm Tr}(\lambda,\t')$ of $g\in G$ such that $g\cdot\lambda\in\t'$ consists of finitely many right cosets mod $T$. A generalised version of Lemma~14.24 (or Lemma~11.9) in \cite{Bo} is also valid if we consider the action of $H$ on an irreducible variety $X$ and replace $\mf{m}$ by an irreducible $M$-stable subvariety of $X$ which has dimension $\dim(X)-\dim(H/M)$. This generalised lemma, applied to the $G$-variety $\g^*$, the subgroup $T$ of $G$ and the subvariety $\t'$ of $\g^*$, gives us that the union of the conjugates of $\t'$ contains an open dense subset of $\g^*$. So, after replacing $T$ by a conjugate, we may assume that $f\vert_{\t'}\ne0$. Now we use the fact that the adjoint action of $\t$ on $\t'\subseteq\g^*$ and $K[\t']$ is trivial and deduce, as in the proof of Lemma~\ref{lem.semiinv}, that $\t\cdot f=0$.\\
(3).\ If $G=\SL_n$, then $G$ inherits a filtration from $K[\Mat_n]$ with associated graded $K[\Mat_n]/(\det)$ which is an integral domain. This filtration is clearly stable under the conjugation actions of $\SL_n$ and $\spl_n$. The fact that every semi-invariant of $\spl_n$ in $K[\SL_n]$ is an invariant follows from Lemma~\ref{lem.semiinv}.

\end{proof}

The fact that $K[\GL_n]^{\gl_n}=K[\gl_n]^{\gl_n}[\det^{-1}]$ is a UFD follows immediately from Theorem~\ref{thm.ufd} or from the fact that this holds for $K[\gl_n]^{\gl_n}$ (see \cite{PrT} Lemma~2).

\begin{remgl}
Assume we have connected normal subgroups $G_1,\ldots,G_k$ of $G$ such that $(G_i,G_j)=1$ for $i\ne j$ and such that the multiplication homomorphism $\prod_{i=1}^kG_i\to G$ is surjective. Denote the Lie algebra of $G_i$ by $\g_i$. For $h\in G$ and $i\in\{1,\ldots,k\}$ define the algebra homomorphism $ev_i^h:K[G]\to K[G_i]$ by $ev_i^h(f)(g)=f(hg)$ for $g\in G_i$. If $h\in\prod_{j\ne i}G_j$, then $ev_i^h$ is equivariant for the conjugation actions of $G_i$ and $\g_i$. Now assume that Proposition~\ref{prop.stableideal} is proved for the $G_i$. Then we can prove it for $G$ as follows. As in the proof we assume that $f\ne0$. Choose a maximal torus $T$ on which $f$ does not vanish and put $\t=\Lie(T)$. As in the proof we get that $\t\cdot f=0$. Now let $i\in\{1,\ldots,k\}$ and $h\in\prod_{j\ne i}G_j$. Then the ideal $K[G_i]ev_i^h(f)$ is $\g_i$-stable, so for $x\in\g_i$ we have that $ev_i^h(x\cdot f)=x\cdot ev_i^h(f)=0$. Since this holds for any $h\in\prod_{j\ne i}G_j$, we must have that $\g_i\cdot f=0$. The conclusion now follows from the fact that $\g=\t+\sum_{i=1}^k\g_i$.

In view of Lemma~\ref{lem.regsselt} this shows that Proposition~\ref{prop.stableideal} and Theorem~\ref{thm.ufd} are also valid if we replace the condition "assume that $\g$ has regular semi-simple elements" by "If $p=2$, assume that the rootsystem of $G$ has no irreducible component of type $C_n$, $n\ge2$", since in Proposition~\ref{prop.application}(3) we have proved those results also for $\SL_2$ in characteristic $2$. Note that this condition holds if $p$ is good for $G$. It also follows that the condition on the existence of regular semi-simple elements in Proposition~\ref{prop.stableideal} and Theorem~\ref{thm.ufd} can be omitted if Proposition~\ref{prop.stableideal} can be proved for $\Sp_{2n}$, $n\ge2$, in characteristic $2$. In \cite{T2} we will obtain this result as a consequence of results on the symplectic ideal.
\end{remgl}

Finally we mention a "rationality result". Contrary to the corresponding statement for $K[\spl_n]^{\spl_n}$ in \cite{PrT} (Thm~1), there are no assumptions on $p$.

\begin{propgl}
$K[\SL_n]^{\spl_n}$ has a rational field of fractions.
\end{propgl}

\begin{proof}
The proof is very similar to (but easier than) that of Theorem~3 in \cite{T1}. By the theorem of \cite{D}, $K[\SL_n]^{\spl_n}$ is generated by $K[\SL_n]^{\SL_n}$ and $K[\SL_n]^p$. Since $K[\SL_n]^{\spl_n}$ is a finitely generated module over $K[\SL_n]^p$ we have that the transcendence degree of the field of fractions of $K[\SL_n]^{\spl_n}$ is the same as that of $K[\SL_n]^p$ which equals $n^2-1$. For $i\in\{1,\ldots,n\}$ define $s_i\in\mb{C}[\Mat_n]$ by $s_i(A)=\tr(\wedge^iA)$, where $\wedge^iA$ denotes the $i$-th exterior power of $A$ and $\tr$ denotes the trace. For $f\in K[\Mat_n]$ denote the restriction of $f$ to $\SL_n$ by $f'$. Then $s'_1,\ldots,s'_{n-1}$ generate $K[\SL_n]^{\SL_n}$ (see e.g. \cite{T1} 1.3). Denote the standard coordinate functionals on $\Mat_n$ by $\xi_{ij}$, $1\le i,j\le n$ and put $\zeta_{ij}=\xi_{ij}^p$. We consider the $s_i$ as polynomials in the $\xi_{ij}$. Then $K[\SL_n]^{\spl_n}$ is generated by the $n^2+n-1$ elements $s'_1,\ldots,s'_{n-1}$, $(\zeta_{ij})_{ij}$ and we have the $n$ relations
$$s_i((\zeta_{ij})_{ij})=s'_i, i=1,\ldots,n-1\text{ and }\det((\zeta_{ij})_{ij})=1.$$
Now we proceed as in the proof of Theorem~3 in \cite{T1}: we observe that these equations are linear in $\zeta_{1n},\ldots,\zeta_{nn}$ and eliminate those $n$ elements from the $n^2+n-1$ generators. So the field of fractions of $K[\SL_n]^{\spl_n}$ is generated by $n^2-1$ elements and therefore purely transcendental.
\end{proof}

The fact that $K[\GL_n]^{\gl_n}$ has a rational field of fractions is an immediate consequence of the fact that this holds for $K[\gl_n]^{\gl_n}$ (see \cite{PrT} Theorem~1, $\gl_n\cong\gl_n^*$ as $\GL_n$-modules).
\medskip

\noindent{\it Acknowledgement}. I would like to thank J.\ C.\ Jantzen for helpful discussions and the Department of Mathematics of the University of Aarhus for its hospitality. This research was partially funded by the EPSRC Grant EP/C542150/1.

\end{document}